\documentclass[11pt]{article}
\usepackage{graphicx}
\usepackage{amssymb}
\usepackage{amsmath,amsthm,amscd,amssymb}
\usepackage{latexsym}
\numberwithin{equation}{section}

\theoremstyle{plain}
\newtheorem{theorem}{Theorem}[section]
\newtheorem{lemma}[theorem]{Lemma}
\newtheorem{corollary}[theorem]{Corollary}

\theoremstyle{definition}
\newtheorem{definition}[theorem]{Definition}

\theoremstyle{remark}

\newtheorem{case[theorem]}{Case}

\def \R{{\mathbb R}}
\def\rr{{\mathbb R}}

\def\supp{\hbox{supp\,}}

\title{Geometric incidence theorems via Fourier analysis}
\author{Alex Iosevich, Hadi Jorati, and Izabella \L aba}
\date{February 4, 2007}

\begin{document}

\maketitle

\begin{abstract} We show that every non-trivial Sobolev bound for
generalized Radon transforms which average functions over families of
curves and
surfaces yields an incidence theorem for suitably regular
discrete sets of points and curves or surfaces in Euclidean space.
This mechanism
allows us to deduce geometric results not readily accessible by
combinatorial methods.
\end{abstract}

%\tableofcontents

\section{Introduction}

An incidence between a set of points $P$ and a set of curves or
surfaces ${\cal G}$ in an Euclidean
space $\rr^d$ is a pair $(p,\Gamma)$, where $p\in P$, $\Gamma\in{\cal
G}$, and $p$
lies on $\Gamma$.   Combinatorial geometers have long been interested
in bounds on
the number of incidences between point sets and families of curves or
a specified type.
Such bounds, apart from their intrinsic interest, have found
applications to a variety of
other combinatorial problems, see e.g. \cite{APST04}, \cite{CEGSW90},
\cite{ST01}, \cite{ST05}; a comprehensive
survey of the area is given in \cite{PS04}.
A prototype result here is the classical Szemer\'edi-Trotter incidence
theorem (\cite{ST83}),
which states that the total number of incidences $I$ between $N$
points and $M$ straight
lines in the plane obeys
\begin{equation}\label{st}
I \lesssim N+M+{(NM)}^{\frac{2}{3}}.
\end{equation}
There are explicit examples showing that this bound is sharp.
Here and below, $c,C$ are constants; $X \lesssim Y$ means that $X \leq
CY$ for some $C$,
and $|\cdot|$ denotes the cardinality of a finite set, the Euclidean
norm of a vector in
$\rr^d$, or the absolute value of a complex number, depending on the
context. In addition,
$X\lessapprox Y$ means $X\leq C_\epsilon N^\epsilon Y$ for any $\epsilon>0$.

Sz\'ekely \cite{Sz97} observed that the Szemer\'edi-Trotter theorem
could be proved using only
very limited geometrical information about straight lines, namely that a
line is
uniquely determined by a pair of points and that two lines can
intersect at only one point.
This allowed him to obtain an extension of the theorem to more general
curves which
satisfy similar intersection axioms.

A more substantial generalization was given by Pach and Sharir
\cite{PS98}, who proved
a version of the Szemer\'edi-Trotter theorem for {\em pseudolines
with $k$ degrees of freedom}. The latter are defined to be a family of
curves in $\R^2$
such that: (a) at most $O(1)$ curves can pass through any given $k$
points, (b) any pair of
curves intersects in at most $O(1)$ points. Under these assumptions,
the number of
incidences between $N$ curves and $M$ points in $\R^2$ is bounded by
$$I\lesssim M^{\frac{k}{2k-1}}N^{\frac{2k-2}{2k-1}}+M+N.$$

%
%More precisely, Sz\'ekely proved the following.
%\begin{theorem} \label{st} (\cite{ST83}, \cite{Sz97}). Let $({\cal
%G},{P})$ be an arrangement
%\footnote{By the arrangement we further mean an embedding, or drawing
%of the curves and points in the plane.}of $M$ points and $N$ curves in
%${\mathbb R}^2$. Suppose that no more than $t$ curves pass through any
%pair of points of ${P}$ and that any two curves of ${\cal G}$
%intersect at no more than $s$ points of ${P}$. Then the total number
%of incidences obeys the bound
%%$$ I \lesssim {(ts)}^{\frac{1}{3}} {(MN)}^{\frac{2}{3}}+tM+N.$$
%\end{theorem}

Formulating and proving similar incidence theorems in higher
dimensions turns out to be
surprisingly difficult, and despite considerable amount of interesting
work in this direction,
comprehensive and sharp estimates are difficult to come by.
To illustrate some difficulties involved in the problem of counting
incidences between points and
codimension-one manifolds -- henceforth surfaces -- consider a family
of $M$ 2-dimensional planes
in $\rr^3$ all containing a fixed line $L$, and a set of $N$ points
all of which lies on $L$.  Then the
number of incidences is $MN$, since every point is incident to every
line.  Thus there can be no
non-trivial incidence bounds for points and planes in $\rr^3$, or for
points and
$k$-dimensional affine subspaces of $\rr^d$ with $k\geq 2$, without
additional assumptions.
Similarly, consider a set of $N$ points in ${\mathbb R}^4$
supported on the unit circle in the $(x_1,x_2)$ coordinate plane, and
a set of $M$ spheres
of radius $\sqrt{2}$ whose centers are supported on the unit circle in
the $(x_3,x_4)$
coordinate plane. Then each point lies on each sphere, and hence the
number of incidences
is again $MN$.   Even in the ``translation invariant" setting,
where the set of ``surfaces" is obtained by translating a single
hypersurface by elements of a
fixed point set, it is still possible to construct similar examples,
for instance an infinite
one-dimensional family of translates of the paraboloid in
${\mathbb R}^3$ given by the equation $x_3=x_1^2+x_2^2$, all
intersecting in a fixed parabola
in a plane parallel to the $x_3$ direction.

These simple examples suggest that incidence theorems for surfaces of
dimension $k>2$
in $\rr^d$ must involve additional geometric assumptions, so as to
preclude lower-dimensional
obstructions such as those just described.   Such results have indeed
been obtained by
combinatorial methods. This was done, for instance in \cite{CEGSW90},
\cite{SV04}, \cite{ST05}
for spheres, in \cite{ET05}, \cite{ST05} for $k$-dimensional affine
subpaces, and in \cite{LS05}
for more algebraic affine surfaces.

The main goal of this paper is to develop a Fourier analytic approach
to the study of incidence
problems for a reasonably general class of families of manifolds of
codimension 1. Our results
depend on certain regularity bounds of operators which average
functions over families of
smooth curves and surfaces in Euclidean space. The properties of such
operators have been
studied for many years in the context of harmonic analysis and PDEs.
See, for example, \cite{St93}
and the references therein for a thorough description of the subject
area. We shall see that
$L^2$-Sobolev  bounds for averaging operators can be converted to
incidence theorems for corresponding surfaces and ``homogeneous" point
sets.   We will also
give specific examples of curves and surfaces to which our incidence
bounds apply.  The main
requirement is that the said regularity improving bounds should hold,
and that has been
proved by harmonic analysts in a wide variety of settings, typically
under assumptions
involving smoothness and curvature.
(By contrast, combinatorial methods  tend to ignore properties of this
type altogether, relying
instead on geometric assumptions similar to those stated above for
pseudolines.)
The homogeneity assumption on the point set,
formulated rigorously in Definition \ref{def-homogeneous}, means that
the point set is uniformly distributed throughout a fixed cube and
eliminates the lower-dimensional
obstructions described earlier.  This condition originated in analysis
literature and was
used in the context of incidence bounds e.g. in \cite{SV04}, \cite{ST05},
\cite{LS05}.

We mention that another interesting type of geometric obstructions
occurs in counting incidences
between points sets and manifolds of higher codimension.   For
example, the Szemer\'edi-Trotter
bound (\ref{st}) extends trivially to point-line arrangements in
higher dimensions, as can be seen
by projecting the arrangement on a fixed plane and observing that each
incidence in the
original arrangement corresponds to an incidence in the projected one.
However, the converse
of the last statement is false, and one expects that the number of
incidences in a truly
higher-dimensional arrangement should be substantially lower.  This is
indeed confirmed
e.g. in \cite{SV04}.  The challenge is thus to find a good notion of a
``higher-dimensional"
arrangement under which sharp bounds, or at least bounds better than
those in two dimensions,
can be obtained.
See, for example, \cite{AKS02}, \cite{LS05}, \cite{SW02}, \cite{Sh94}
for the case of incidences between
points and lines or curves in higher dimensions. Connections between
this type of problems
and certain questions in analysis is described and thoroughly
referenced in \cite{W99};
see also \cite{Sch03}, which puts some recent work in harmonic
analysis involving combinatorics
of circles and spheres in a broader perspective.
We do not address this class of problems here, but we do setup a
framework which we hope
to apply to such questions in the future.

%%%%%%%%%%%%%%%%%%%%%%%%%%%%%%%%%%

\section{The main result}

%%%%%%%%%%%%%%%%%%%%%%%%%%%%%%%%%%%

Our results concern upper bounds on the number of incidences between
``homogeneous"
point sets (see below) and families of surfaces for which certain
analytic bounds are known
to hold.  A model case that the reader should keep in mind is that of
a finite and homogeneous
point set $A\subset [0,1]^d$ and the family of spheres $\{S_a\}_{a\in
A}$, where
$S_a=\{x\in {\mathbb R}^d:\ |x-a|=.1\}$. We will obtain an upper bound
on the number of pairs
$\{(a,a')\in A\times A: a'\in S_a\}$ (equivalently, $\{(a,a')\in
A\times A: |a-a'|=.1\}$).  While this
particular case has already been investigated by combinatorial methods
and we are not able
to improve the known bounds here, our results apply to other classes
of surfaces as well, for
example spheres with varying radii (depending on $a$) or translates of
smooth convex
hypersurfaces satisfying appropriate analytic conditions.  We do not
know of any combinatorial
methods that would work in the latter case -- simple extensions of the
known methods are not
sufficient, basically because it is difficult to control intersections
of convex bodies in
dimensions 3 and higher.  There is some overlap between the results
presented here and
those of {\L}aba and Solymosi \cite{LS05}, where combinatorial assumptions
on
the surfaces are made instead of analytic assumptions.

We first describe the point sets that we work with.
\begin{definition}\label{def-homogeneous}
Let $C_0,c_0$ be positive constants with $0<c_0<C_0$.
We say that a set $A \subset {[0,1]}^d$, $d \ge 2$, is {\em
$(C_0,c_0,k_0)$-homogeneous}
if every cube of sidelength $c_0|A|^{-1/d}$ contains at most $k_0$
points of $A$, and
if every cube of sidelength
$C_0|A|^{-1/d}$ contains at least one point of $A$.
\end{definition}

We recall that an infinite point set $A \subset {\mathbb R}^d$ is
called a {\em Delone}
(a.k.a. {\em Delauney} or {\em well-distributed}) set if  there exist
$0<c_A<C_A$ such
that every cube of sidelength $c_A$ contains at most $k_0$ points of $A$,
and
every cube of sidelength $C_A$
contains at least one point of $A$. Thus the rescaled and truncated sets
$$A_t=[0,1]^d\cap\{ta:a\in A\},\ $$ are $(C_A,c_A,k_0)$-homogeneous,
with the same
constants for all $t>1$.
Delone sets have been recently studied in connection with higher-dimensional
incidence theorems and with the Erd\H{o}s and Falconer distance
problems, see e.g.
\cite{IH}, \cite{IL2}, \cite{SV04}, \cite{I04}, \cite{ST05}, \cite{IR06}.

The Fourier transform of a Schwartz class function $f$ is defined by the
formula
$$\hat{f}(\xi)=\int_{{\mathbb R}^d} e^{-2 \pi i x \cdot \xi} f(x)dx.$$
In what follows, the homogeneous Sobolev space $L^2_{\gamma}({\mathbb R}^d)$
is the closure of the Schwartz class ${\cal S}({\mathbb R}^d)$ in the norm
$${||u||}^2_{L^2_{\gamma}({\mathbb R}^d)}=
\int_{{\mathbb R}^d} {|\xi|}^{2\gamma} {|\widehat{u}(\xi)|}^2 d\xi. $$

\renewcommand \theenumi{\roman{enumi}}

\begin{definition}\label{gamma-def}
Let ${\cal G}=\{\Gamma_x\}_{x \in (0,1)^d}$ be a continuous family of
smooth manifolds in  ${\mathbb R}^{d}$.

\begin{enumerate}
\item We say that ${\cal G}$ is {\em regular} if there is a smooth
function $\Phi:\rr^{2d}_{x,y}\to\rr$
such that $\Gamma_x=\{y:\ \Phi(x,y)=0\}$ for all $x\in[0,1]^d$, and that
$|\nabla_x\Phi(x,y)|>\epsilon_0$, $|\nabla_y\Phi(x,y)|>\epsilon_0$ for some
$\epsilon_0>0$ and all $(x,y)$ with $\Phi(x,y)=0$.
\item
We say that ${\cal G}$ is $\gamma$-regular if the averaging operator
$T$, given by
$$ Tf(x)=\int_{\Gamma_x} f(y) d\sigma_x(y), $$ where $\sigma_x$ is the
Lebesgue
measure\footnote{
We do not normalize $\sigma_x$, so that the measure of $\Gamma_x$ is
equal to $\liminf_{\delta\to
0}|\Gamma_x^\delta|$. } on $\Gamma_x$, obeys the following Sobolev estimate:
\begin{equation}\label{sobolev}
{||Tf||}_{L^2_{\gamma}({\mathbb R}^d)} \lesssim {\|f\|}_{L^2({\mathbb
R}^d)}.
\end{equation}
%\item \item We say that ${\cal G}$ is locally $\gamma$-regular if the
%averaging operator $T$, given above satisfies
%%$$ {||Tf_j||}_{L^2({\mathbb R}^d)} \lesssim 2^{-j \gamma}
%{||f_j||}_{L^2({\mathbb R}^d)},$$ where $f_j$ is the Littlewood-Paley
%piece of $f$ defined by $\widehat{f}_j(\xi)=\widehat{f}(\xi)
%\phi_0(2^{-j}|\xi|)$, where $\phi_0$ is a smooth cut-off function
%supported in $[1/2,4]$ and equal to $1$ in $[1,2]$.
%\item
%Let $1\leq p<q<\infty$.  We say that $\mathcal{G}$
%is strongly $(p,q)$-regular if the operator $T$ defined above
%obeys the following $L^p$-improving estimate:
%\begin{equation}\label{pq-improve}
%\| Tf \|_{L^q(\mathbb{R}^d)} \lesssim \| f \|_{L^p(\mathbb{R}^d)}.
%\end{equation}
\item
We say that ${\cal G}$ is strongly $\gamma$-regular if the following
holds for all $t\in(-c_R,c_R)$,
where $c_R$ is a small fixed positive constant.
Let $\Gamma_{x,t}=\{y:\Phi(x,y)=t\}.$ Then the families
${\cal G}_t=\{\Gamma_{x,t}\}_{x\in(0,1)^d}$ are $\gamma$-regular,
with the constant uniform for all $t\in (-c_R,c_R)$.
%Similarly, we say that ${\cal G}$ is strongly $(p,q)$-regular if
%${\cal G}_t$ are
%$(p,q)$-regular.
%\item Analogously, we say that ${\cal G}$ is locally strongly
%$\gamma$-regular if it is strongly $\gamma$-regular with the condition
%that ${\cal G}_t=\{\Gamma_{x,t}\}_{x\in(0,1)^d}$ are $\gamma$-regular,
%with the constant uniform for all $t\in (-c_R,c_R)$, replaced by the
%same condition with ${\cal G}_t=\{\Gamma_{x,t}\}_{x\in(0,1)^d}$
%locally $\gamma$-regular.

\end{enumerate} \end{definition}

We remark that, since $L^p$ functions are defined only up to sets of
measure zero, the restriction of such a
function to a lower-dimensional submanifold need not always be defined
or measurable. Nonetheless, under
certain natural conditions (examples of which will be given shortly),
$Tf(x)$ is defined for almost all $x$
and obeys (\ref{sobolev}). In fact, for compact surfaces,
(\ref{sobolev}) with $\gamma=0$  trivially holds in
our setup by the Minkowski's integral inequality. Plugging this into
our numerology would result in
a trivial version of Theorem \ref{main} with incidence bound
$N^{2-\frac{1}{d}}$, reflecting the
observation that each $\Gamma_x$ can be incident to at most
$N^{(d-1)/d}$ points due to
dimensionality considerations and the homogeneity assumption.
It is easy to construct counterexamples showing that the
estimate (\ref{sobolev}) generally cannot hold with
$\gamma>\frac{d-1}{2}$, see e.g. \cite{Ho83}.

We now state our main results.

\begin{theorem} \label{main}
Let $A$ be a $(C_0,c_0,k_0)$-homogeneous set in ${\mathbb R}^d$,
$d \ge 2$, with $|A|=N$.
Suppose that the family ${\cal G}=\{\Gamma_x\}_{x \in (0,1)^d}$ is strongly
$\gamma$-regular for some $\gamma>0$. Let $s\in(d-\gamma,d)$.
Then for all $\delta\leq c N^{-1/s}$, where $c$ is small enough,
\begin{equation}\label{mest} \left|
\left\{(a,a') \in A \times A: a' \in \Gamma_{a}^\delta \right\}
\right|\;\lesssim
N^{2-\frac{1}{s}}.\end{equation}
In particular, we have the incidence bound
\begin{equation}\label{main-2}
\left|
\left\{(a,a') \in A\times A: a' \in  \Gamma_{a} \right\}
\right|\;\lessapprox
N^{2-\frac{1}{d-\gamma}}.\end{equation}

\end{theorem}

Here and through the rest of the paper, the implicit constants in the
$\lesssim$ and $\lessapprox$
symbols depend on the constants $c_0, C_0$ in Definition
\ref{def-homogeneous},
the constants $C_L,c_R, \epsilon_0$ in Definition \ref{gamma-def}, and
the exponent $s$.  However,
they are always independent of $N$.

Note that if $A,B$ are two $(C_0,c_0,k_0)$-homogeneous sets, then $A\cup B$
is
$(C_0,c_0,2k_0)$-homogeneous.  Therefore our results apply just as
well (possibly
with different constants) to quantities of the form
$$|\left\{(a,b) \in A\times B: b \in  \Gamma_{a} \right\}|,$$
where $A,B$ are two different $(C_0,c_0,k_0)$-homogeneous sets.  This
comment
applies to Theorem \ref{main} as well as to all corollaries in
the sequel.

\section{Applications}\label{sec-incidences}

In this section, we discuss applications of the abstract Theorem \ref{main}
to specific incidence problems.
We continue to assume that $A\subset[0,1]^d$ is a
$(C_0,c_0,k_0)$-homogeneous
set with $|A|=N$.

We first discuss a group of results obtained by combining Theorem
\ref{main} and classical theory
of Fourier integral operators. We need a definition.

\begin{definition}\label{monge-ampere}
Suppose that ${\cal G}=\{\Gamma_x\}_{x\in{\mathbb R}^d}$ is a family
of  $d-1$-dimensional hypersurfaces in
${\mathbb R}^d$ defined by $\Gamma_x=\{y: \Phi(x,y)=0\}$, where
$\Phi:{\mathbb R}^d\times{\mathbb
R}^d\to{\mathbb R}$ is a smooth function.  We will say that ${\cal G}$
satisfies the rotational curvature
condition of Phong and Stein (see \cite{PS86} and \cite{So93}, Theorem
6.2.1 and Corollary 6.2.3) if the
Monge-Ampere determinant
\begin{equation} {\cal M}(\Phi)\;=\;\left( \begin{matrix} 0 &
\frac{\partial\Phi}{\partial x_1}&\ldots & \frac{\partial\Phi}{\partial x_d}
\\
\frac{\partial\Phi}{\partial y_1}& \frac{\partial^2\Phi}{\partial
x_1\partial y_1}&\ldots&
\frac{\partial^2\Phi}{\partial x_d\partial y_1}
\\
 \vdots&\vdots&\ddots &\vdots
\\
\frac{\partial\Phi}{\partial y_d}& \frac{\partial^2\Phi}{\partial
x_1\partial y_d}&\ldots&
\frac{\partial^2\Phi}{\partial x_d\partial y_d}
\end{matrix} \right),
\label{mad}\end{equation}
restricted to the set where $\Phi(x,y)=0$, does not vanish.
\end{definition}

\begin{corollary}\label{cor-ma1}
Suppose that ${\cal G}=\{\Gamma_x\}_{x\in{\mathbb R}^d}$, where
$\Gamma_x=\{y: \Phi(x,y)=0\}$,  satisfies the
Phong-Stein condition (\ref{mad}) in Definition \ref{monge-ampere}.
Let $s\in(\frac{d+1}{2},d)$.
Then for any $\delta<cN^{-1/s}$, with $c$ small enough,
\begin{equation}\label{ma-e0}
\left|\left\{(a,a') \in A\times A: a' \in \Gamma_a^{\delta} \right\}\right|
\;\lesssim N^{2-\frac{1}{s}}.
\end{equation}
In particular,
\begin{equation}\label{ma-e1}
\left|\left\{(a,a') \in A\times A: a' \in  \Gamma_a \right\}\right|
\;\lessapprox\; N^{\frac{2d}{d+1}}.
\end{equation}
\end{corollary}

Corollary \ref{cor-ma1} follows from Theorem \ref{main} in view of the
main result of \cite{PS86} which says
that under the assumptions of Corollary \ref{cor-ma1}, the estimate
(\ref{sobolev}) holds with
$\gamma=\frac{d-1}{2}$.  Note also that the Phong-Stein condition
implies that $|\nabla_x\Phi|$
and $|\nabla_y\Phi|$ are bounded from below away from 0 on $\Gamma_x$.
The constants hidden in the $\lesssim$ and $\lessapprox$ symbols
will clearly depend on the function $\Phi$.

While the Monge-Ampere condition in Corollary \ref{cor-ma1} is often
interpreted as a curvature
condition, it also allows for certain families of ``flat" hypersurfaces.
For example, it is satisfied by the family of hyperplanes
$\Gamma_x=\{y \in {(0,1)}^d: x \cdot y=1\}$.
The corresponding estimate in Corollary \ref{cor-ma1} yields a special
case of the
Szemer\'edi-Trotter theorem in two-dimensions and a non-trivial incidence
bound
in higher dimension.

As a special case of Corollary \ref{cor-ma1} we obtain the following.

\begin{corollary}\label{cor-ma2}
Let $d \ge 2$, and let $r(x)$
be a smooth function $[0,1]^d\to(0,\infty)$ such that $|\nabla
r(x)|\leq c<1.$  Let
$\Gamma_x=\{y:|x-y|=r(x)\}$.  Then the conclusions of Corollary
\ref{cor-ma1} hold,
in particular we have
\begin{equation}\label{ma-e2}
\left|\left\{(a,a') \in A\times A:  |a-a'|= r(a)\right\}\right|\;
\lessapprox\;
N^{\frac{2d}{d+1}}.
\end{equation}
\end{corollary}

To prove the corollary, it suffices to verify that the Phong-Stein
condition holds for
$\Phi(x,y)=|x-y|^2-r(x)^2$.  We do this in Section
\ref{proof-MA}.

If $d=2$ and $r(x)\equiv r_0$ is fixed, Corollary \ref{cor-ma2} says
in particular that the number
of pairs $a,a'\in A$ such that $|a-a'|=r_0$ (i.e. the number of
incidences between the $N$ points
of $A$ and $N$ circles of radius $r_0$ and
centered at those points)  is $\lessapprox N^{4/3}$.  This is a
partial result on the
``unit distance" conjecture of Erd\H{o}s, which asserts that
if $P\subset{\mathbb R}^2$ is a set of cardinality $N$ and $r_0>0$ is
fixed, then
\begin{equation}\label{unit-dist}
\# \{(a,b) \in P \times P: |a-b|=r_0\} \lesssim N \sqrt{\log(N)}.
\end{equation}
The best known partial result to date follows from
the Szemer\'edi-Trotter theorem and yields $N^{4/3}$ on the right side
of (\ref{unit-dist})
\cite{SST}. Our result matches this up to the endpoint.
On the one hand, our point set is somewhat special and we do lose the
endpoint. 
On the other hand, our theorem applies also to circles $\{y:|a-y|=r(a)\}$
of varying radii, as well as to thin annuli
$\{y: r(a)\leq |a-y|\leq r(a)+\delta\}$ with $\delta=cN^{-2/3}$.  These extensions
do not appear to follow
from the known combinatorial theorems (Szemer\'edi-Trotter, Pach-Sharir)
in any straightforward way.

In three dimensions, the best known estimate on the number of unit distances
is
$O(N^{\frac{3}{2}+\epsilon})$ for any $\epsilon>0$ \cite{CEGSW90}.
Our result matches
this estimate in the case of homogeneous sets, except that \cite{CEGSW90}
gives
an explicit form of the endpoint $o(N^{\epsilon})$ factor.
(The best known lower bound is $N^{4/3}\log\log N$.)  In dimensions $d\geq
4$,
the number of unit distances in a general point set of cardinality $N$
can be of the order
$N^2$, as
demonstrated by the example in Section 1.  Nonetheless, Corollary
\ref{cor-ma2} still
yields a non-trivial bound for homogeneous sets in higher dimensions.

We note that the condition $|\nabla r(x)|<c<1$ has an appealing geometrical interpretation:
it is a slightly strengthened quantitative version of the statement that
that no sphere $\Gamma_x$ is entirely contained within another sphere $\Gamma_y$.

We further note that our theorem yields the same conclusions for the
unit distance
problem if the Euclidean norm is
replaced by a non-isotropic norm
in which the unit sphere is a convex set with sufficiently smooth boundary
and
everywhere non-vanishing curvature.
It is well known that combinatorial methods run into
difficulties for this type of problems, especially in higher
dimensions.  The known
incidence bounds for spheres use very specific geometric information,
for example
that two spheres always intersect along a circle; on the other hand,
the intersection curves of more general convex surfaces can be almost
impossible to control.
The related paper \cite{LS05} gives an incidence bound for algebraic
surfaces in
$\rr^3$, with exponent depending on the algebraic degree of the
surfaces.  The bounds
in Corollaries \ref{cor-ma1}, \ref{cor-ma2} improve on that of
\cite{LS05} for algebraic
surfaces of high enough degree.

In Corollary \ref{cor-ma2}, we assumed that the radius $r(x)$ of the
sphere centered at $x$
obeyed $|\nabla r(x)|<c<1$. If we assume instead that $|\nabla
r(x)|>C>1$ for all $x$,
then it turns out that a weaker estimate can be proved
in dimensions $d\geq 3$.  Again, our estimate is in fact more general
and applies
just as well to translated and dilated copies of a fixed curved
hypersurface.  We first give the analytic statement of the result.

\begin{corollary}\label{cor3}
Let $d \ge 3$. Let ${\cal G}=\{\Gamma_x\}_{x \in [0,1]^d}$, where
$\Gamma_x=\{y: \Phi(x,y)=0\}$ and 
$\Phi:{\mathbb R}^d\times{\mathbb
R}^d\to{\mathbb R}$ is a smooth function.
Assume that:

(i) $|\nabla_x\Phi(x,y)|>\epsilon_0$, $|\nabla_y\Phi(x,y)|>\epsilon_0$ for some
$\epsilon_0>0$ and all $(x,y)$ with $\Phi(x,y)=0$,

(ii) there is a smooth function $r(x)$ with values in $[a,b]$ for some $0<a<b<\infty$
such that
$\Phi(x,y)=f(\frac{y-x}{r(x)})$ for some smooth function $f:\rr^d\to\rr$, so that $\Gamma_x
=x+r(x)\Gamma$ for a fixed $\Gamma:=\Gamma_0$, 

(iii) $\Gamma$ (defined above) is a smooth closed
hypersurface with everywhere non-vanishing Gaussian curvature.

Let $s\in(\frac{d+2}{2},d)$.
Then for any $\delta<cN^{-1/s}$, with $c$ small enough,
\begin{equation}\label{ma-e3}
\left|\left\{(a,a') \in A\times A: a' \in \Gamma_a^{\delta} \right\}\right|
\;\lesssim N^{2-\frac{1}{s}}.
\end{equation}
In particular,
\begin{equation}\label{ma-e4}
\left|\left\{(a,a') \in A\times A: a' \in  \Gamma_a \right\}\right|
\;\lessapprox\; N^{\frac{2d+2}{d+2}}.
\end{equation}
\end{corollary}

We now give a variant of the above corollary which is slightly weaker, but easier to
apply to combinatorial problems where the defining function $\Phi$ is not given
explicitly.  

\begin{corollary}\label{cor4}
Let $d \ge 3$. Let $K$ be a convex body in $\rr^d$ whose interior contains 0 and
such that $\Gamma=\partial K$  is a smooth closed
hypersurface with everywhere non-vanishing Gaussian curvature.
Let $\Gamma_x=x+r(x)\Gamma$, where $r(x)$
is a smooth function with values in $[a,b]$ for some $0<a<b<\infty$.
Assume that either
\begin{equation}\label{tt.e1}
|\nabla r(x)|>C_0>M \hbox{ for all }x\in[0,1]^d
\end{equation}
or 
\begin{equation}\label{tt.e2}
|\nabla r(x)|<c_0<m \hbox{ for all }x\in[0,1]^d,
\end{equation}
where $m=\min_{x\in\Gamma} |x|$ and $M=\max_{x\in\Gamma}|x|$.
Let $s\in(\frac{d+2}{2},d)$.
Then for any $\delta<cN^{-1/s}$, where $c$ is small enough,
the estimates (\ref{ma-e3}) and (\ref{ma-e4}) hold.
\end{corollary}

In particular, the assumptions (therefore the conclusion) of the corollary hold in the special 
case when $\Gamma_x=\{y:|x-y|=r(x)\}$ and $r(x)$ is a smooth function with 
$|\nabla r(x)|>C_0>1$ for all
$x$, since then $\Gamma$ is the unit sphere and $m=M=1$. (The case 
$|\nabla r(x)|<c_0<1$ is already covered by Corollary \ref{cor-ma2}.)

If $d=2$, the allowed range of $s$ in the first estimate (\ref{ma-e3})
is empty.  In this
case, the estimate (\ref{ma-e4}) is the trivial incidence bound
which follows from
the dimensionality of $\Gamma_x$ and homogeneity of $A$.

\section{Proof of Theorem \ref{main}}\label{proof-main}

The proof is based on a conversion mechanism, first developed in
\cite{IH}, \cite{IL2}
in the context of the Falconer distance set problem.
Let $A$ and $\{\Gamma_x\}_{x \in (0,1)^d}$ be as in the statement of
the theorem.
Let also $\delta=cN^{-1/s}$, where $c>0$ is small enough and
\begin{equation}\label{e-gamma-s}
d-\gamma<s<d.
\end{equation}
We will prove that
\begin{equation}\label{e-main}
\left| \left\{(a,a') \in A \times A: a' \in  \Gamma_{a}^\delta
\right\}\right| \;\lesssim\;
N^{2-\frac{1}{s}}.
\end{equation}
This clearly implies the theorem, since the inequality clearly remains
valid if the $\delta$
on the left side of (\ref{e-main}) is replaced by a smaller number.

We define
\begin{equation}
f(x)=N^{-1}\delta^{-d} \sum_{a\in A}
\phi\left(\frac{x-a}{\delta}\right),
\label{dens}\end{equation}
where $\phi:{\mathbb
R}^d\to[0,\infty)$ is a smooth function such that $\phi(x)\equiv 1$
for $|x|\leq 1$ and $\phi(x)\equiv 0$ for
$|x|\geq 2$.  Then $f$ is supported on a $\delta$-neighbourhood of $[0,1]^d$
and
$\int_{\R^d}f\approx 1$.  Let $d\mu=fdx$. We also let
$$E=\bigcup_{a\in A}\{x:|x-a|\leq \delta\}.$$

For  each pair $(a,a') \in A\times A$, let
$$B_{a,a'}:=\left\{(x,y): \left|x-a\right|\leq c\delta,
\ \left|y-a'\right| \leq c\delta \right\}.$$
We will assume that the constant $c<1$ is small enough so that the
sets $B_{a,a'}$ are pairwise disjoint.
We also have $B_{a,a'}\subset E \times E$.

We further observe that if $c$ is sufficiently small (which we will
assume henceforth),
and if $(x,y) \in B_{a,a'}$ and $a' \in \Gamma_{a}^{c\delta}$, then $
y\in \Gamma_{x}^{\delta}$.
This follows from the uniform continuity of the defining function $\Phi$.
Hence, with $c$ as above,
\begin{equation}\label{count}
N^{-2}\left|\left\{(a,a') \in A\times A: a' \in \Gamma_a^{c \delta}
\right\}\right|
\;\lesssim\;  \mu\times \mu\,\{(x,y) \in E\times E: y \in
\Gamma_x^{\delta} \}.\end{equation}

Let
\begin{equation}\label{fops}\begin{array}{lll}
T^{\delta}g(x)&=&\int_{\Gamma_x^{\delta}}
g(y)dy, \\ \hfill \\
T_tg(x)&=&\int_{\Gamma_{x,t}} g(y) d\sigma_{x,t}(y),
\end{array}\end{equation} where $\Gamma_{x,t}$ are as in
Definition \ref{gamma-def} and $d\sigma_{x,t}$ is the surface measure
on $\Gamma_{x,t}$.
Using a change of variables, and invoking the regularity of $\Phi$
again, we estimate
\begin{equation}\label{zz.e1}
\begin{array}{lll} \mu\times
\mu\,\{(x,y) \in E\times E: y \in \Gamma_x^{\delta} \}&=& \langle
T^{\delta}f(x),f(x)\rangle \\ \hfill \\
&\lesssim& \int_0^{C\delta} \langle T_tf(x),f(x)\rangle dt. \end{array}
\end{equation}
We now use (\ref{sobolev}) to bound the last integrand uniformly in $t$.

Fix Schwartz class functions $\eta_0(\xi)$ supported in $|\xi|\leq 4$
and $\eta(\xi)$
supported in the spherical shell $1<|\xi|<4$  such that the quantities
$\eta_0(\xi)$,
$\eta_j(\xi)=\eta(2^{-j}\xi),\,j\geq1$ form a partition of unity.

Write $f=\sum_{j=0}^\infty f_j$, where
$\widehat{f_j}(\xi)=\widehat{f}(\xi)\eta_j(\xi)$. Then
\begin{equation}
\langle T_tf,f\rangle =\sum_{j,k}\langle T_tf_j,f_k\rangle
=\sum_{|j-k|\leq K}+\sum_{|j-k|>K},
\label{ps}\end{equation}
where $K$ is a large enough constant.  We will estimate the two sums
separately,
starting with the second one.  Note that this is very easy in the
translation-invariant
case when $\Gamma_x=\Gamma+x$ for a fixed $\Gamma$ (i.e. $\Phi(x,y)$
depends only
on $x-y$), since then, by Plancherel's theorem,
$$\langle T_tf_j,f_k\rangle =\langle
\widehat{T_tf_j},\widehat{f_k}\rangle =\langle
\widehat{f_j*\sigma_{0,t}},\widehat{f_k}\rangle$$
$$=\langle \widehat{f_j}\widehat{\sigma_{0,t}},\widehat{f_k}\rangle
=0\hbox{ if }|j-k|>1$$
and the second sum vanishes if we let $K=1$.
In the general case, we need the following lemma.

\begin{lemma} \label{integrationbyparts} Assume that $K>0$ is large enough.
Then for any $M$ there exists $C_M>0$ such that for all $k,j$ with $|j-k|>K$
$$ \langle Tf_j, f_k\rangle  \leq C_M2^{-M\max(j,k)}.$$
\end{lemma}

We defer the proof of the lemma until the next section.  Applying the lemma
with $M=1$, we estimate the second sum in (\ref{ps}):
\begin{equation}
\sum_{|j-k| > K}|\langle T_tf_j,f_k\rangle|
\lesssim \sum_{j=0}^\infty \sum_{k=j+K+1}^\infty 2^{-k}
+\sum_{k=0}^\infty \sum_{j=k+K+1}^\infty 2^{-j}
\lesssim 1.
\label{ps1}\end{equation}

Turning to the first sum in (\ref{ps}), we write
$$\sum_{|j-k|\leq K}|\langle T_tf_j,f_k\rangle|
=\sum_{r=-K}^K \sum_{j=0}^\infty|\langle T_tf_j,f_{j+r}\rangle|,
$$
where we put $f_k\equiv 0$ for $k<0$.
By (\ref{sobolev}) and the Cauchy-Schwartz inequality
we have
\begin{equation}\label{cse}|\langle T_tf_j,f_{j+r}\rangle|
\leq \|T_tf_j\|_2\|f_{j+r}\|_2
\leq
2^{-j\gamma}\|f_j\|_2\|f_{j+r}\|_2.
\end{equation}
We claim that
\begin{equation}\label{l2norm}
\|f_j\|_2^2\lesssim 2^{j(d-s)}.
\end{equation}
Assuming (\ref{l2norm}) for the moment, we conclude that
\begin{equation}\label{ps2}
\sum_{|j-k|\leq K}|\langle T_tf_j,f_k\rangle|
\lesssim\sum_{r=-K}^K \sum_{j=0}^\infty
2^{-j\gamma}2^{j(d-s)/2}2^{(j+r)(d-s)/2}
\lesssim \sum_{j}2^{-j\gamma+j(d-s)}\lesssim 1,
\end{equation}
provided that (\ref{e-gamma-s}) holds.
Combining (\ref{ps1}) and (\ref{ps2}), we see that
$$|\langle T_tf,f\rangle|\lesssim 1.$$

Finally, we plug this into (\ref{zz.e1}) and get
$$\mu\times \mu\,\{(x,y) \in E\times E: y \in \Gamma_x^{\delta}
\}\lesssim\delta
\lesssim N^{-1/s}.$$
The conclusion (\ref{e-main}) follows from this and (\ref{count}).

It remains to prove (\ref{l2norm}).  By the Fourier support
localization of $f_j$ and
Plancherel's theorem, we have
$$\|f_j\|_2^2=\|\widehat{f_j}\|_2^2
\;\lesssim \;2^{j(d-s)}\int|\widehat{f}(\xi)|^2|\xi|^{-d+s}d\xi.$$  On
the other hand, the integral in the
right-hand-side is the $s$-energy of $\mu$, i.e.
$$\int|\widehat{f}(\xi)|^2|\xi|^{-d+s}d\xi
=c_{d,s}\int\int |x-y|^{-s}d\mu(x)d\mu(y):=I_s,$$
where $c_{d,s}$ is an explicit constant depending only on $d$
and $s$ (see \cite{W03}).
Thus it suffices to prove that
\begin{equation}
I_s\lesssim 1.
\end{equation}
Indeed, we write
$$I_s=N^{-2}\delta^{-2d}\sum_{a,a'\in A}\int\int |x-y|^{-s}
\phi \left(\frac{x-a}{\delta}\right)\phi \left(\frac{y-a'}{\delta}\right)$$
$$=N^{-2}\delta^{-2d}\left(\sum_{a,a'\in A:|a-a'|\leq 4\delta}
+\sum_{j=2}^\infty \sum_{a,a'\in A:2^j\delta<|a-a'|\leq 2^{j+1}\delta}
\right):= N^{-2}\delta^{-2d}(S_1+S_2).$$
We start with the first term.  If $\delta\lesssim N^{-1/s}$ with
$0<s<d$, then for
large $N$ we can only have $|a-a'|<4\delta$, $a,a'\in A$, if $a=a'$.  Hence
$$
S_1=\sum_{a}\int\int |x-y|^{-s}
\phi \left(\frac{x-a}{\delta}\right)\phi \left(\frac{y-a}{\delta}\right)
\lesssim N \int\int |u|^{-s}\phi \left(\frac{u+y}{\delta}\right)\phi
\left(\frac{u}{\delta}\right)$$
$$\lesssim N \delta^d \int_{|u|\leq 2\delta}\int |u|^{-s}ds
\lesssim N\delta^{2d-s}.$$
In the second term, if $2^j\delta<|a-a'|\leq 2^{j+1}\delta$, then
$2^{j-1}\delta<|x-y|\leq 2^{j+2}\delta$ on the support of the integrand.
For each $a\in A$, there are about $N\cdot(2^j\delta)^d$ points $a'\in A$
with
$2^j\delta<|a-a'|\leq 2^{j+1}\delta$.
Hence
$$S_2\lesssim \sum_{j=2}^\infty N^2(2^j\delta)^d
(2^j\delta)^{-s}\delta^{2d}$$
$$\lesssim \sum_{j=2}^\infty N^2\delta^{3d-s}2^{j(d-s)}\lesssim
N^2\delta^{3d-s},$$
where we again used that $0<s<d$.  Combining the estimates on $S_1$ and
$S_2$,
we get
$$I_s\lesssim N^{-2}\delta^{-2d}(N\delta^{2d-s}+N^2\delta^{3d-s})
=N^{-1}\delta^{-s}+\delta^{d-s}\lesssim 1,$$
as claimed.
This proves
(\ref{l2norm}) and completes the proof of Theorem \ref{main}.

%
%%%%%%%%%%%%%%%%%%%%%%%%%%%%%%%%%%%%%%%%%%%%%%%

\section{Proof of Lemma \ref{integrationbyparts}}

To simplify the notation, we will only prove the lemma with $T_t$
replaced by the operator
$T$ as in Definition \ref{gamma-def}.  It will be clear from the proof
that the same estimates
hold for $T_t$ for $|t|\leq C\delta$, with constants uniform in $t$.

We write
$$Tf(x)=\int_{\{y: \Phi(x,y)=0\}} f(y) \psi(x,y)d\sigma_x(y),$$
where $\psi$ is smooth and compactly supported and $d\sigma_x(y)$ is
the surface
measure on $\Gamma_x=\{y: \Phi(x,y)=0\}$. Using the
$\gamma$-regularity assumption,
and shrinking the support of $\psi$ if necessary, we may
assume that $C^{-1}\leq|\nabla_y\Phi(x,y)|\leq C$ for some $C>0$ and all
$x,y\in
\supp\psi$.
By the definition of Lebesgue measure, we
can approximate $T$ by $T_n$ as $n\to\infty$:
$$
T_nf(x)=n \int f(y) \psi(x,y) \psi_0(n\Phi(x,y))dy,$$
where $\psi_0$ is a smooth cut-off function supported in $[-1,1]$.
Hence is suffices to
prove our estimate for $T_n$ with constants uniform in $n$.  Applying
Fourier inversion
twice, to $\psi_0$ and then to $f$, we write
$$ T_nf(x)=\int \int e^{2 \pi i t \Phi(x,y)} f(y) \psi(x,y)
\widehat{\psi}_0(n^{-1}t) dtdy.$$
$$=\int \int \int e^{2\pi i\mu\cdot y} e^{2 \pi i t \Phi(x,y)}
\widehat{f}(\mu) \psi(x,y)
\widehat{\psi}_0(n^{-1}t) dtdyd\mu.$$It follows that
$$\widehat{T_nf}(\xi)=
\int \int \int \int  e^{2\pi i\mu\cdot y}e^{2 \pi i t
\Phi(x,y)}\widehat{f}(\mu) \psi(x,y) \widehat{\psi}_0(n^{-1}t)
e^{-2\pi ix\cdot\xi} dtdydxd\mu.$$
Using Plancherel's theorem, we write
$\langle T_nf_j,f_k\rangle$ as
\begin{equation} \label{parad} \begin{array}{ll}
&\langle T_nf_j,f_k\rangle=\langle \widehat{T_nf_j},
\widehat{f}_k\rangle\\ \hfill \\
&=\int \int \int \int \int e^{2 \pi i t \Phi(x,y)} e^{-2 \pi i x \cdot
\xi} e^{2 \pi i y \cdot \mu} \widehat{f}_j(\mu)
\widehat{f}_k(\xi) \psi(x,y) \widehat{\psi}_0(n^{-1}t) dt dy dx d\mu
d\xi\\ \hfill \\
&=\int\int\int
\widehat{f_j}(\mu)\widehat{f_k}(\xi)I_{jk}(\xi,\mu,t)\widehat{\psi}_0(n^{-1}
t)
,
\end{array}\end{equation}
where
$$ I_{jk}(\xi,\mu,t)=\phi_0(2^{-j}|\mu|) \phi_0(2^{-k}|\xi|) I(\xi,\mu,t)$$,
\begin{equation}\label{ii-1}
I(\xi,\mu,t)=\int \int e^{2 \pi i t \Phi(x,y)} e^{-2 \pi i x \cdot \xi}
e^{2 \pi i y \cdot \mu} \psi(x,y)dxdy
\end{equation}
and $\phi_0$ is a fixed function in $C_0^\infty(\rr^d)$ equal to 1 on
$\{1\leq |x|\leq 10\}$ and
vanishing on $\{|x|\leq\frac{1}{2}\}$.

We claim that if $|j-k|>K$, where $K$ is a large enough constant, then
for any $M$ there is a constant $C_M$ such that
\begin{equation}\label{ijk-decay}
|I_{jk}(\xi,\mu,t)|\leq C_M 2^{-M\max(j,k)},
\end{equation}
uniformly in $\xi,\mu,t$.
Assuming this, the proof of the lemma is completed as follows.  By
(\ref{parad}) and
(\ref{ijk-decay}), we have
$$
|\langle T_nf_j,f_k\rangle|\lesssim
C_M2^{-M\max(j,k)}\|\widehat{f_j}\|_1\|\widehat{f_k}\|_1.
$$
By the Cauchy-Schwartz inequality and  (\ref{l2norm}),
$$
\|\widehat{f_j}\|_1\leq \|\widehat{f_j}\|_2|\supp f_j|^{1/2}\lesssim
2^{jd/2}2^{jd/2}\lesssim 2^{jd},
$$
so that
$$
|\langle T_nf_j,f_k\rangle|\lesssim C_M2^{-M\max(j,k)}2^{jd}2^{kd}
\lesssim C_M 2^{-(M-2d)\max(j,k)}.
$$
Relabelling the constants, we get the conclusion of the lemma.

It remains to prove (\ref{ijk-decay}).
The idea is that $I$ is an oscillatory integral with critical points given
by
$$ t \nabla_x \Phi(x,y)=\xi; \ \ t \nabla_y \Phi(x,y)=-\mu.$$
If $|j-k|\geq K$, where $K$ is large enough depending on $c_0,C_0$,
then at least one of the cut-off functions $\phi_0(2^{-j}|\mu|)$,
$ \phi_0(2^{-k}|\xi|)$
is supported away from the critical points, hence we can estimate $I_{jk}$
using
the easy part of the stationary phase method which only involves integration
by
parts.

The details are as follows.
Recall that we are assuming that
\begin{equation}\label{gradients}
c_0\leq |\nabla_x \Phi(x,y)|\leq C_0, \
c_0\leq |\nabla_y \Phi(x,y)|\leq C_0,
\end{equation}
for some positive constants $c_0,C_0$.
We have
\begin{equation} \label{byparts-1}
\frac{1}{2 \pi i \left(t \frac{\partial \Phi}{\partial x_m}+\xi_m \right)}
\frac{\partial}{\partial x_m} e^{2 \pi i (t \Phi(x,y)+x \cdot
\xi)}=e^{2 \pi i (t \Phi(x,y)+x \cdot \xi)}, \end{equation} and,
similarly,
\begin{equation} \label{byparts-2} \frac{1}{2 \pi i \left(t
\frac{\partial \Phi}{\partial y_m}-\mu_m \right)}
\frac{\partial}{\partial y_m} e^{2 \pi i (t \Phi(x,y)-y \cdot
\mu)}=e^{2 \pi i (t \Phi(x,y)-y \cdot \mu)}.
\end{equation}
Using (\ref{byparts-1}) in (\ref{ii-1})
and integrating $M$ times by parts in $x_m$, we see that for any $M$ there
is a constant $C'_M$ (depending on $\Phi$ and the cut-off functions
$\psi,\phi_0$ but not
on $\xi,\mu,t$ or $j,k$) such that for $m=1,\dots,d$
$$
| I(\xi,\mu,t)|\leq C'_M \max_{x,y}
\left|t \frac{\partial \Phi}{\partial x_m}+\xi_m \right|^{-M}.
$$
If $|t|\leq (2C_0)^{-1}|\xi|$, where $C_0$ is as in (\ref{gradients}),
then $|t \nabla\Phi|\leq \frac{|\xi|}{2}$, hence
there is at least one $m$ such that
\begin{equation}\label{ii-2}
\left|t \frac{\partial \Phi}{\partial x_m}+\xi_m \right|
\geq d^{-1}\Big||t\nabla\Phi|-|\xi|\Big|
\geq \max\Big(\frac{|\xi|}{2d},
\frac{|t|}{C_0d}\Big).
\end{equation}
Similarly, if $|t|\geq 2c_0^{-1}|\xi|$, then $|t \nabla\Phi|\geq 2|\xi|$,
hence
\begin{equation}\label{ii-2b}
\left|t \frac{\partial \Phi}{\partial x_m}+\xi_m \right|\geq
\max\Big(\frac{|\xi|}{d},
\frac{c_0|t|}{2d}\Big)
\end{equation}
for at least one $m$.  It follows that in both regions,
\begin{equation}\label{ii-3}
| I_{jk}(\xi,\mu,t)|\leq C''_M (\max(|t|,|\xi|))^{-M}.
\end{equation}
Similarly, if either $|t|\leq (2C_0)^{-1}|\mu|$ or $|t|\geq 2c_0^{-1}|\mu|$,
we have for any $M$
\begin{equation}\label{ii-4}
| I_{jk}(\xi,\mu,t)|\leq C''_M (\max(|t|,|\mu|))^{-M}.
\end{equation}
Thus one of (\ref{ii-3}), (\ref{ii-4}) must hold unless we have both
$(2C_0)^{-1}|\xi|\leq |t|\leq  2c_0^{-1}|\xi|$
and $(2C_0)^{-1}|\mu|\leq |t|\leq  2c_0^{-1}|\mu|$.
But then
\begin{equation}\label{ii-5}
\frac{1}{4}c_0C_0^{-1}|\mu|\leq |\xi|\leq 4c_0^{-1}C_0|\mu|.
\end{equation}
Choose $K$ large enough so that $K\geq 100 c_0^{-1}C_0$.  Then (\ref{ii-5})
fails on the support of $\phi_0(2^{-j}|\mu|) \phi_0(2^{-k}|\xi|)$ whenever
$|j-k|\geq K$, so that at least one of (\ref{ii-3}), (\ref{ii-4}) hold.
Suppose now that $j-k>K$, so that $|\xi|\leq|\mu|$ on $\supp I_{jk}$.
If (\ref{ii-4}) holds, (\ref{ijk-decay}) follows immediately.  If on the
other hand (\ref{ii-4}) fails, we must in particular have $|t|\geq
(2C_0)^{-1}|\mu|$.
Plugging this into (\ref{ii-3}), we get (\ref{ijk-decay}) again.

%%%%%%%%%%%%%%%%%%%%%%%%%%%%%%%%%%%%%%%%%%%%%%%

\section{Proof of Corollary \ref{cor-ma1}}\label{proof-MA}

It suffices to show that ${\cal M}(\Phi)$, where
\begin{equation}
\Phi(x,y)=(x_1-y_1)^2+\dots+(x_d-y_d)^2-r(x)^2,
\label{3.1}\end{equation} has determinant bounded away from zero
on the set $\{(x,y):\ \Phi(x,y)=0\}$. We abbreviate
$r_j=\frac{\partial r}{\partial x_j}$. Then
\begin{equation}\label{3.2}\begin{array}{lll}
\det[{\cal M}(\Phi)]&=&\left| \begin{matrix} 0 &
2(y_1-x_1)&2(y_2-x_2)&\ldots & 2(y_d-x_d)
\\
2(x_1-y_1)-2rr_1& -2 &0 &\ldots& 0
\\
2(x_2-y_2)-2rr_2& 0 &-2 &\ldots& 0
\\
 \vdots&\vdots&\vdots &\ddots &\vdots
\\
2(x_d-y_d)-2rr_d& 0 &0 &\ldots& -2
\end{matrix} \right|\\ \hfill \\
&=&(-1)^d 2^{d+1} D_d,\end{array}
\end{equation}
where
\begin{equation}\label{3.3}
\hspace{-15mm}D_d\;\;=\;\;\left| \begin{matrix} 0 &
x_1-y_1&x_2-y_2&\ldots & x_d-y_d
\\
x_1-y_1-rr_1& 1 &0 &\ldots& 0
\\
x_2-y_2-rr_2& 0 &1 &\ldots& 0
\\
 \vdots&\vdots&\vdots &\ddots &\vdots
\\
x_d-y_d-rr_d& 0 &0 &\ldots& 1
\end{matrix} \right|.
\end{equation}
Expanding in the last row, we get
\begin{equation}\label{3.4}
D_d=(-1)^d(x_d-y_d-rr_d) \left| \begin{matrix} x_1-y_1&x_2-y_2&\ldots
&x_{d-1}-y_{d-1}& x_d-y_d
\\
1 &0 &\ldots&0& 0
\\
0 &1 &\ldots& 0&0
\\
 \vdots&\vdots&\ddots &\vdots &\vdots
\\
0 &0 &\ldots& 1&0
\end{matrix} \right|+D_{d-1}.
\end{equation}
Expanding the remaining determinant in the last column yields
\begin{equation}\label{3.5}
D_{d} =-(x_d-y_d)^2+rr_d(x_d-y_d)+D_{d-1}.
\end{equation}
It follows by induction that
\begin{equation}\label{3.6}
D_d=-|x-y|^2+r(x)(x-y)\cdot \nabla r(x).
\end{equation}
Therefore, on the set where $\Phi(x,y)=0$, we have
\begin{equation}\label{3.7}
D_d=-r(x)^2+r(x)(x-y)\cdot \nabla r(x).
\end{equation}
We rewrite this as
\begin{equation}\label{3.8}
D_d=r(x)^2(r(x)^{-1}(x-y)\cdot \nabla r(x)-1).
\end{equation}
Note that $r(x)^{-1}(x-y)$ is a unit vector on the set where $\Phi(x,y)=0$,
and
that we are assuming that
$|\nabla r(x)|<c<1$. Hence $|r(x)^{-1}(x-y)\cdot \nabla r(x)-1|\geq 1-c>0$,
as claimed. This completes the proof.

%%%%%%%%%%%%%%%%%%%%%%%%%%%%%%%%%%%%%%%%

\section{Proof of Corollary \ref{cor3}}\label{proof-cor34}

Corollary \ref{cor3} follows from the proof of Theorem \ref{main} and the estimate, implicit in
\cite{Gr81} (see also \cite{St93}), that if
$$ {\cal A}_{\Gamma}f(x)=\sup_{1<t<2} \left| \int_{\Gamma} f(x-ty)
d\sigma(y) \right|, $$ where $\Gamma$ is a smooth hypersurface with
non-vanishing Gaussian curvature and
$d\sigma$ is the Lebesgue measure on $\Gamma$, then
\begin{equation}\label{dd.e1}
 \|\mathcal{A}_{\Gamma}f_j\|_{L^2(\mathbb{R}^d)} \lesssim
2^{-\frac{j(d-2)}{2}}\|f_j\|_{L^2(\mathbb{R}^d)},
\end{equation}
where $f_j$
is the Littlewood-Paley piece of $f$ defined as in Section \ref{proof-main} (after 
(\ref{zz.e1})).  
Since the averaging operator $T$ associated with $\mathcal{G}$ is dominated by 
the maximal operator $\mathcal{A}_\Gamma$, it follows that 
$ \|Tf_j\|_{2} \lesssim 2^{-\frac{j(d-2)}{2}}\|f_j\|_{2}$,  We now plug this 
directly into (\ref{cse}) in the proof of Theorem \ref{main} (instead of using the
assumption (\ref{sobolev})).  Since this is the only place where (\ref{sobolev}) was
used, and all remaining assumptions of the theorem are satisfied, we
obtain the corollary.

We now prove (\ref{dd.e1}). We have
$$ \mathcal{A}_{\Gamma}f(x):= \sup_{1<t<2} |\mathcal{A}_t f(x)| \qquad
\textrm{where} \qquad A_t f(x) := \int f(x-ty) d\sigma(y).$$ 
It suffices to show that
\begin{equation} \label{thekeyestimate} \int \sup_t |A_t f_j(x)|^2 dx
\leq C \|f_j\|_{2}^{2}.(2^j)^{-(d-2)}.\end{equation}
In order to prove the claim, first observe that by the Fourier
inversion formula we have
\begin{equation} \label{fourierform} A_t f_j (x) = \int_{|\xi| \sim
2^j} e^{2 \pi i x \cdot \xi}
\widehat{f_j}(\xi)\widehat{\sigma}(t\xi)d\xi. \end{equation}
By the classical method of stationary phase (see e.g. \cite{So93} and
comments in \cite{IS97}),
\begin{equation} \label{classicaldecay} |\widehat{\sigma}(\xi)| \leq
C|\xi|^{-\frac{d-1}{2}}, \end{equation} and, similarly,
\begin{equation} \label{classicaldecay2} | \xi \cdot \nabla
\widehat{\sigma}(\xi)| \leq C|\xi| \cdot
{|\xi|}^{-\frac{d-1}{2}}. \end{equation}

We need the following basic lemma which is proved using the
fundamental theorem of calculus and the Cauchy-Schwartz inequality.
See, for example, \cite{So93}, Chapter 2.
\begin{lemma} \label{sobolevbyhand} Let $F$ be a continuously
differentiable function. Then
$$ \sup_{t \in [1,2]} {|F(t)|}^2 \leq {|F(1)|}^2+2{\left( \int_1^2
{|F(t)|}^2 dt \right)}^{\frac{1}{2}} \cdot {\left( \int_1^2
{|F'(t)|}^2 dt \right)}^{\frac{1}{2}}. $$
\end{lemma}

We now apply  Lemma (\ref{sobolevbyhand}) to $F(t)=A_tf_j(x)$ and use
Cauchy-Schwartz to see that
$$ \int \sup_{t \in [1,2]} {|A_tf_j(x)|}^2 dx \leq \int {|A_1f_j(x)|}^2 dx$$
$$+ {\left( \int \int_1^2 {|A_tf_j(x)|}^2 dxdt \right)}^{\frac{1}{2}}
\cdot {\left( \int \int_1^2
{\left|\frac{d}{dt} A_tf_j(x)\right|}^2 dxdt \right)}^{\frac{1}{2}}$$
$$=I+II \cdot III.$$
By Plancherel and (\ref{classicaldecay}),
$$ I \leq C{||f_j||}_2^2 \cdot 2^{-j(d-1)}.$$
Applying Plancherel once again,
$$ II \leq C {||f_j||}_2 \cdot 2^{-j \frac{d-1}{2}}.$$
Since $\frac{d}{dt} \widehat{\sigma}(t \xi)=\xi \cdot
\widehat{\sigma}(t \xi),$ Plancherel and (\ref{classicaldecay2}) imply
that
$$ III \leq C {||f_j||}_2 \cdot 2^j \cdot 2^{-j \frac{d-1}{2}}.$$
Combining these estimates we see that
$$ |A_t f_j(x)|^2 \leq C \|f_j\|_{2}^{2} \cdot 2^j \cdot
(2^j)^{-(d-1)}, $$ where $C$ does not depend on $t$. It follows that
$$ \left( \sup_{t \in [1,2]} |A_t f_j (x)|^2 dx \right)^{\frac{1}{2}}
\leq C \|f_j\|_2.(2^j)^{-\frac{d-2}{2}}, $$
which is what we wanted to prove.

%%%%%%%%%%%%%%%%%%%%%%%%%%%%%%%%%%%%%
\section{Proof of Corollary \ref{cor4}}\label{proof-cor4}

Recall that the Minkowski functional of $K$ is the unique function $\phi_K:\rr^d\to[0,\infty)$
which is 
homogeneous degree 1 and equal to 1 on $\Gamma$.  Let $f(x)$ be a smooth
function on $\rr^d$, non-vanishing except on $\Gamma$, such that 
$f(x)=1-\phi_K(x)$ for $\frac{am}{2}<|x|<2bM$.
We would like to apply Corollary \ref{cor3} with $\Phi(x,y)=f(\frac{y-x}{r(x)})$.  It is
clear from the definition that all assumptions are satisfied, except for (i) which 
needs to be verified.  By explicit computation, we have for all $x,y$ with
$|\Phi(x,y)|<\epsilon$ small enough
\begin{equation}\label{tt.e3}
\nabla_y\Phi(x,y)=-\frac{1}{r(x)}\nabla f\Big(\frac{y-x}{r(x)}\Big),
\end{equation}
\begin{equation}\label{tt.e4}
\nabla_x\Phi(x,y)=\frac{1}{r(x)}\Big(
\nabla f\Big(\frac{y-x}{r(x)}\Big) +\Big(\frac{y-x}{r(x)}\cdot \nabla f\Big(\frac{y-x}{r(x)}\Big)
\nabla r(x)\Big).
\end{equation}
(\ref{tt.e3}) is clearly bounded away from 0 as required, by the definition of $f$.  Turning 
to (\ref{tt.e4}), we note that since $r$ is bounded from above and below, and since 
$\frac{y-x}{r(x)}\in\Gamma$ whenever $\Phi(x,y)=0$, it suffices to prove that the
quantity
\begin{equation}\label{tt.e5}
\nabla f(u)+(u\cdot\nabla f(u))\nabla r(x)
\end{equation}
has norm bounded away from 0 if $u\in\Gamma$.  Using the definition of $f$, and computing
its directional derivative in the direction of $u$ at $u\in\Gamma$, we get
$u\cdot\nabla f(u)=1$ on $\Gamma$. We also get that for $u\in\Gamma$,
$$
|\nabla f(u)|=\frac{1}{|u|\cos\theta},
$$
where $\theta$ is the angle between $u$ and the outward normal vector to $\Gamma$
at $u$.  But the right side  is also equal to the distance between the origin and the tangent
line to $\Gamma$ at $u$.  By the convexity of $K$, this is minimized when $|u|=m$ and
maximixed when $|u|=M$, and at those points we have $\cos\theta=1$.  Thus
$m\leq |\nabla f(u)|\leq M$.  It follows that for $u\in\Gamma$,
$$|(\ref{tt.e5})|\geq \big| |\nabla f(u)| -r(x)\big|>\epsilon_0>0$$
if one of (\ref{tt.e1}), (\ref{tt.e2}) holds.  By continuity, a similar estimate (with 
$\epsilon_0$ replaced by $\epsilon_0/2$) holds on a neighbourhood of 
$\Gamma$.  This completes the verification of (i) and hence proves the corollary.

%%%%%%%%%%%%%%%%%%%%%%%%%%%%%%%%%%%%%%

\end{document}